\documentclass{amsart}

\let\hat=\widehat
\let\tilde=\widetilde

\newcommand{\BF}{{\mathbb F}}
\newcommand{\BP}{{\mathbb P}}
\newcommand{\BQ}{{\mathbb Q}}
\newcommand{\BZ}{{\mathbb Z}}

\newcommand{\CC}{{\mathcal C}}
\newcommand{\CM}{{\mathcal M}}

\newcommand{\character}{\operatorname{char}}
\newcommand{\End}{\operatorname{End}}
\newcommand{\Gal}{\operatorname{Gal}}
\newcommand{\Id}{\operatorname{Id}}
\newcommand{\Jac}{\operatorname{Jac}}
\newcommand{\Resultant}{\operatorname{resultant}}

\newcommand{\mybar}[1]{#1\llap{$\overline{\phantom{\rm#1}}$}}

\newcommand{\Kb}{\mybar{K}}

\newcommand{\lra}{\longrightarrow}
\newcommand{\llra}{\relbar\joinrel\longrightarrow}
\newcommand{\lmr}[1]{\mathop{\llra}\limits^{#1}}

\newcommand{\md}[1]{\vbox{\vbox to 4pt{}\vbox{\hbox{
     \Big\downarrow\rlap{$\vcenter{\hbox{$\scriptstyle#1$}}$}}\vfill}}}

\newcommand{\mup}[1]{\vbox{\vbox to 4pt{}\vbox{\hbox{
     \Big\uparrow\rlap{$\vcenter{\hbox{$\scriptstyle#1$}}$}}\vfill}}}

\newtheorem{theorem}{Theorem}

\newtheorem{proposition}[theorem]{Proposition}
\newtheorem{example}[theorem]{Example}

\theoremstyle{definition}
\newtheorem*{defn}{Definition}

\theoremstyle{remark}
\newtheorem*{rem}{Remark}

\begin{document}

\title[Curves with isomorphic Jacobians]
{Infinite families of pairs of curves over $\BQ$
 with isomorphic Jacobians}
 
\author{Everett W.\ Howe}
\address{Center for Communications Research, 
         4320 Westerra Court, 
         San Diego, CA 92121-1967, USA.}
\email{however@alumni.caltech.edu}
\urladdr{http://alumni.caltech.edu/\~{}however/}

\date{12 May 2003}
 


\begin{abstract}
We present three families of pairs of geometrically non-isomorphic
curves whose Jacobians are isomorphic to one another as unpolarized
abelian varieties.  Each family is parametrized by an open subset
of~$\BP^1$.  The first family consists of pairs of genus-$2$ curves
whose equations are given by simple expressions in the parameter; the
curves in this family have reducible Jacobians.  The second family also
consists of pairs of genus-$2$ curves, but generically the curves in
this family have absolutely simple Jacobians.  The third family 
consists of pairs of genus-$3$ curves, one member of each pair being a
hyperelliptic curve and the other a plane quartic.  Examples from these
families show that in general it is impossible to tell from the
Jacobian of a genus-$2$ curve over $\BQ$ whether or not the curve has 
rational points --- or indeed whether or not it has real points.  Our
constructions depend on earlier joint work with Franck Lepr\'evost and
Bjorn Poonen, and on Peter Bending's explicit description of the curves
of genus $2$ whose Jacobians have real multiplication
by~$\BZ[\sqrt{2}]$.
\end{abstract}

\maketitle

\section{Introduction}
\label{S-intro}

Torelli's theorem shows that a curve is completely determined by its
polarized Jacobian variety, but it has been known since the late 1800s
that distinct curves can have isomorphic unpolarized Jacobians.  In
particular, the unpolarized Jacobian of a curve may not reflect all of
the curve's geometric properties.  Proving that a particular property
of curves cannot always be determined from the Jacobian is equivalent
to showing that there exist two curves, one with the given property and
one without, whose Jacobians are isomorphic to one another.  Thus, for
example, the pairs of curves written down in~\cite{Howe:PAMS} show that
one cannot tell whether or not a curve of genus $3$ over the complex
numbers is hyperelliptic simply by looking at its Jacobian.

One would also like to find {\it arithmetic\/} properties of curves
that are not determined by the Jacobian, but from an arithmetic
perspective the heretofore-known explicit examples of distinct curves
with isomorphic Jacobians (catalogued in the introduction 
to~\cite{Howe:PAMS}) are not entirely satisfying.  The primary
complaint is that none of the examples involves curves that can be
defined over $\BQ$; in addition, for any given number field only
finitely many of the examples can be defined over that field.  
Furthermore, all of the explicit examples in characteristic $0$ known
before now involve curves with geometrically reducible Jacobians, and
the arithmetic of such curves differs qualitatively from that of curves
whose Jacobians are irreducible.

In this paper we address these concerns by providing three new explicit
families of pairs of non-isomorphic curves with isomorphic Jacobians.
Each family is parametrized by an open subset of~$\BP^1$, so each
family gives an infinite number of examples over~$\BQ$.  Also, the
Jacobians of the curves in one of the families are typically absolutely
simple.  Using examples from these families, we show that the Jacobian
of a genus-$2$ curve over $\BQ$ does not determine whether or not the
curve has rational points, or indeed whether or not the curve has real
points.  Liu, Lorenzini, and Raynaud~\cite{LLR} have used our results
to show that the Jacobian of a genus-$2$ curve over $\BQ$ does not
determine the number of components on the reduction of a minimal model
of the curve modulo a prime.

Our first family of pairs of curves can be defined over an arbitrary
field $K$ whose characteristic is not~$2$.  If $t$ is an element of $K$
with $t(t+1)(t^2+1)\neq 0$ then the equation 
$$(t + 1) y^2 = (2 x^2 - t) (4 t^2 x^4 + 4(t^2 + t + 1)x^2 + 1)$$
defines a curve of genus $2$ that we will denote~$C(t)$.  Clearly the 
quotient of $C(t)$ by the involution $(x,y)\mapsto (-x,y)$ is an 
elliptic curve, so the Jacobian of $C(t)$ splits over~$K$.
         
\begin{theorem}
\label{T-nonsimple2}
Let $K$ be a field of characteristic not $2$ and suppose $t$ is an 
element of $K$ such that $t(t^2-1)(t^2+1)$ is nonzero.  Then $C(t)$ and
$C(-t)$ are curves of genus~$2$ over $K$ whose Jacobians are isomorphic
over~$K$.  Furthermore, $C(t)$ and $C(-t)$ are geometrically 
non-isomorphic unless $K$ has characteristic $11$ and 
$t^2 \in \{-3, -4\}$.
\end{theorem}

Our next family takes a little more effort to describe.  In order to do
so we must define the {\em Richelot duals\/} of a genus-$2$ curve over
a field $K$ of characteristic not~$2$ 
(see~\cite[Ch.~9]{Cassels-Flynn}, \cite[\S3]{Bending}).  Suppose $C$ is
a genus-$2$ curve over $K$ defined by an equation $\delta y^2 = f$,
where $\delta\in K^*$ and where $f$ is a monic separable polynomial in 
$K[x]$ of degree $6$.  Let $\Kb$ be a separable closure of $K$, and
suppose $f$ can be factored as a product $g_1 g_2 g_3$ of three monic
quadratic polynomials in $\Kb[x]$ that are permuted by $\Gal(\Kb/K)$. 
For each $i$ write $g_i = x^2 - t_i x + n_i$ and suppose the 
determinant
$$ d = \left|
\begin{matrix}
1 & t_1 & n_1 \\
1 & t_2 & n_2 \\
1 & t_3 & n_3
\end{matrix}
\right| ,$$
which is an element of $K$, is nonzero.  Define three new polynomials
by setting
\begin{align*}
h_1 &= g_3 \frac{dg_2}{dx} - g_2 \frac{dg_3}{dx}\\
h_2 &= g_1 \frac{dg_3}{dx} - g_3 \frac{dg_1}{dx}\\
h_3 &= g_2 \frac{dg_1}{dx} - g_1 \frac{dg_2}{dx}.
\end{align*}
Then the product $h_1 h_2 h_3$ is a separable polynomial in $K[x]$ of
degree $5$ or~$6$.

\begin{defn}
The {\em Richelot dual of $C$ associated to the factorization 
$f = g_1 g_2 g_3$} is the genus-$2$ curve $D$ defined by
$d \delta y^2 = h_1 h_2 h_3.$
\end{defn}

\begin{theorem}
\label{T-simple2}
Let $K$ be a field of characteristic not $2$, let $v$ be an element of
$K\setminus\{0,1,4\}$ such that 
$$
(v^2 - v + 4)
(v^2 + v + 2)
(v^2 + 3 v + 4)
(v^3 - 6 v^2 - 7 v - 4)
(v^3 - 4 v^2 + 7 v + 4)
\neq 0,$$
let $w$ be a square root of $v$ in $\Kb$, and define numbers 
$\rho_1,\ldots,\rho_6$ by setting
\begin{align*}
\rho_1 & = \frac{(-2+w)(1+w)}{2w^2}    
       & \qquad \qquad \rho_2 & = \frac{(-2-w)(1-w)}{2w^2}\\
\rho_3 & = \frac{-2(2+w)}{(-2+w)(1+w)} 
       & \qquad \qquad \rho_4 & = \frac{(-2-w)(1-w)}{(-w)(-1-w)} \\
\rho_5 & = \frac{(-2+w)(1+w)}{w(-1+w)} 
       & \qquad \qquad \rho_6 & = \frac{-2(2-w)}{(-2-w)(1-w)}.
\end{align*}
The $\rho_i$ are distinct from one another, so that if we let 
$f = \prod(x-\rho_i)$ then the curve $D$ over $K$ defined by $y^2 = f$
has genus~$2$.  Set 
\begin{align*}
g_1 & = (x-\rho_1)(x-\rho_5)    
    & \qquad \qquad g_1' & = (x-\rho_1)(x-\rho_3)\\
g_2 & = (x-\rho_2)(x-\rho_4)    
    & \qquad \qquad g_2' & = (x-\rho_2)(x-\rho_6)\\
g_3 & = (x-\rho_3)(x-\rho_6)    
    & \qquad \qquad g_3' & = (x-\rho_4)(x-\rho_5).
\end{align*}
The Richelot duals $C$ and $C'$ of $D$ with respect to the 
factorizations $f = g_1 g_2 g_3$ and $f = g_1' g_2' g_3'$ exist, and
their Jacobians become isomorphic to one another over 
$K(\sqrt{v(v-4)}\,)$.   The curves $C$ and $C'$ are geometrically
non-isomorphic unless one of the following conditions holds{\/\rm:}
\begin{itemize}
\item[(a)] $\character K = 3$ and 
           $v^{10} - v^8 + v^7 - v^6 - v^5 + v + 1 = 0${\/\rm;}
\item[(b)] $\character K = 19$ and $v + 1 = 0${\/\rm;}
\item[(c)] $\character K = 89$ and $v + 36 = 0${\/\rm;}
\item[(d)] $\character K = 1033$ and $v + 508 = 0$.
\end{itemize}
Furthermore, if $K$ has characteristic $0$ and if $v$ is not an 
algebraic number, then the Jacobians of $C$ and $C'$ are absolutely 
simple.
\end{theorem}

In fact, when $K$ has characteristic $0$ it is very easy to find 
{\it algebraic\/} numbers $v$ in $K$ for which the Jacobians in 
Theorem~\ref{T-simple2} are absolutely simple.  For example, suppose 
$R$ is a subring of $K$ for which there is a homomorphism $\varphi$
to an extension of~$\BF_{13}$.  We show in the proof of
Theorem~\ref{T-simple2} that in this case the Jacobians of $D$ and 
$D'$ are geometrically irreducible whenever $v$ lies in
$\varphi^{-1}(2)$ or~$\varphi^{-1}(6)$. 

Theorem~\ref{T-simple2} gives a $1$-parameter family of pairs of 
non-isomorphic curves with isomorphic Jacobians.  In fact, we shall
see that there is a family of such pairs of curves parametrized by an
elliptic surface; over $\BQ$, this surface has positive rank.

Our third family of pairs of curves with isomorphic Jacobians is again
easy to write down.  Suppose $K$ is a field of characteristic not $2$
and suppose $t$ is an element of $K$ with 
$t(t+1)(t^2+1)(t^2+t+1)\neq 0$.  Let $H(t)$ be the genus-$3$
hyperelliptic curve defined by the homogeneous equations
\begin{align}
W^2 Z^2 &= - \frac{ (t^2+1) }{t(t+1)(t^2+t+1)}X^4
           - \frac{4(t^2+1) }{ (t+1)(t^2+t+1)}Y^4
           + \frac{        1}{t              }Z^4 \label{E-H1}  \\
     0  &= - X^2 + 2 t Y^2 + (t+1) Z^2            \label{E-H2}
\end{align}
and let $Q(t)$ be the plane quartic
\begin{multline}
\label{E-Q}
X^4 + 4t^2 Y^4 + (t+1)^2 Z^4  + (8t^2 + 4t + 8)X^2Y^2 \\ 
       - (4t^2 + 2t + 2)X^2Z^2 + (4t^2 + 4t +8)Y^2Z^2 = 0.
\end{multline}

\begin{theorem}
\label{T-nonsimple3}
Let $K$ be a field of characteristic not $2$ and let $t$ be an element
of $K$ such that $t(t+1)(t^2+1)(t^2+t+1)\neq 0$.  Then the Jacobians of
the two genus-$3$ curves $H(t)$ and $Q(t)$ are isomorphic to one 
another over~$K$.
\end{theorem}

In Section~\ref{S-curves} we mention some simple facts about abelian
surfaces with two non-isomorphic principal polarizations and we show 
how Richelot isogenies can in principle be used to produce such
surfaces from an abelian surface that has nontrivial automorphisms.
In Section~\ref{S-nonsimple} we prove Theorem~\ref{T-nonsimple2}. 
In Section~\ref{S-RM} we review a result of Bending that shows how to
obtain genus-$2$ curves over a given field $K$ whose Jacobians have 
real multiplication by~$\sqrt{2}$ over $K$, and we show how to adapt
Bending's result to obtain curves over $K$ with real multiplication 
by~$\sqrt{2}$ over a quadratic extension of~$K$.  In 
Section~\ref{S-Galois} we give some Galois restrictions on our 
generalization of Bending's construction that ensure that the curves
we construct have two rational Richelot isogenies to curves with 
isomorphic Jacobians. In Section~\ref{S-application} we show that there
is a positive-rank elliptic surface whose points give rise to pairs of
genus-$2$ curves with isomorphic Jacobians, and we prove 
Theorem~\ref{T-simple2}.  In Section~\ref{S-genus3} we prove 
Theorem~\ref{T-nonsimple3}.  Finally, in Section~\ref{S-examples} we 
provide some explicit examples of curves over $\BQ$ produced by our 
theorems, and we show that the Jacobian of a curve over $\BQ$ does
not determine whether or not the curve has rational points, or even
whether or not it has real points.

We relied heavily on the computer algebra system Magma~\cite{Magma}
while working on this paper.  Some of our Magma routines are available
on the web:  to find them, start at

\smallskip
\noindent
{\tt http://alumni.caltech.edu/\~{}however/biblio.html}
\smallskip

\noindent
and follow the links related to this paper.

\section{Abelian surfaces with non-isomorphic polarizations}
\label{S-curves}

Weil~\cite{Weil} showed that an abelian surface with an indecomposable
principal polarization is a Jacobian, so one of our goals in this paper
is to write down abelian surfaces with two non-isomorphic principal 
polarizations.  In this section we will make a few observations about
such surfaces.

Suppose $B$ is an abelian surface with two principal polarizations
$\mu$ and~$\mu'$, which we view as isogenies from $B$ to its dual 
variety~$\hat{B}$.  The polarized varieties $(B,\mu)$ and $(B,\mu')$
are isomorphic to one another if and only if there is an automorphism
$\beta$ of $B$ such that $\mu' = \hat{\beta} \mu \beta$, where 
$\hat{\beta}$ is the dual of~$\beta$.  We would like to write down an
abelian surface $B$ with two non-isomorphic principal polarizations 
$\mu$ and~$\mu'$, so we would like to avoid the existence of such an 
automorphism~$\beta$.  We will accomplish this by obtaining $\mu'$ 
from $\mu$ through the use of an automorphism of a surface 
{\it isogenous\/} to~$B$.  Our main tool is the following well-known 
construction:

Suppose $(A,\lambda)$ is a principally-polarized abelian surface over
a field~$K$, suppose $n$ is a positive integer, and suppose $G$ is a 
rank-$n^2$ subgroupscheme of the $n$-torsion $A[n]$ of $A$ that is 
isotropic with respect to the $\lambda$-Weil pairing on~$A[n]$.  Let 
$B$ be the quotient abelian surface $A/G$ and let $\varphi:A\to B$ be
the natural map.  Then there is a unique principal polarization $\mu$ 
of $B$ that makes the following diagram commute:
$$\begin{matrix}
      A      & \lmr{n\lambda} &      \hat{A}        \\
\md{\varphi} &                & \mup{\hat{\varphi}} \\
      B      &   \lmr{\mu}    &      \hat{B}
\end{matrix}$$

Now suppose that $A$ has an automorphism $\alpha$ such that 
$G':=\alpha(G)$ is also an isotropic subgroup of $A[n]$, and let 
$(B',\mu')$ be the principally polarized abelian surface obtained from 
$G'$ as above.

\begin{proposition}
\label{P-basic}
The automorphism $\alpha$ of $A$ provides an isomorphism $B\to B'$.
If we identify $B'$ with $B$ via this automorphism, then 
$\mu' = \hat{\beta} \mu \beta$, where $\beta$ is the image of
$\alpha^{-1}$ in $(\End B)\otimes \BQ$.
\end{proposition}

\begin{proof}
We have a commutative diagram with exact rows:
$$\begin{matrix}
0 &\lra &      G     &\lra &      A     &\lmr{\varphi} &B  &\lra & 0\\
  &     &\md{\alpha} &     &\md{\alpha} &              &   &     &  \\
0 &\lra &      G'    &\lra &      A     &\lmr{\varphi'}&B' &\lra & 0,
\end{matrix}$$
where $\varphi$ and $\varphi'$ are the natural maps from $A$ to $B$ 
and~$B'$, respectively.  Completing the diagram, we find an isomorphism
$B\to B'$.  This proves the first statement of the theorem.  The second
statement follows by an easy diagram chase.
\end{proof}

This proposition leaves us some hope, because the $\beta$ in the 
proposition will not be an element of $\End B$ if $G\neq G'$.   Also,
if we consider the case $n=2$ and if the principally-polarized surface 
$(A,\lambda)$ is given to us explicitly as either a Jacobian or a
product of polarized elliptic curves, then the theory of the Richelot
isogeny will allow us to write down $(B,\mu)$ and $(B,\mu')$ explicitly
as Jacobians, as we explain below.  Thus, we would like to explicitly
write down abelian surfaces $A$ with non-trivial automorphisms.  In 
later sections we will consider two families of such explicitly-given 
surfaces:  products of isogenous elliptic curves, and Jacobians with 
real multiplication by~$\sqrt{2}$.

We close this section with a comment about Richelot duals and maximal
isotropic subgroups.  Suppose $C$ is a genus-$2$ curve defined by an 
equation $\delta y^2 = f$ and suppose $D$ is the Richelot dual of $C$ 
corresponding to a factorization $f = g_1 g_2 g_3$.  Then there is an
isogeny from the Jacobian of $C$ to the Jacobian of $D$ whose kernel 
is the order-$4$ subgroup $G$ of $\Jac C$ containing the classes of the
divisors $(a_i,0) - (b_i,0)$, where $a_i$ and $b_i$ are the roots of 
$g_i$ in~$\Kb$ (see~\cite[Ch.~9]{Cassels-Flynn}).  The subgroup $G$ is
a maximal isotropic subgroup of the $2$-torsion of $\Jac C$ under the
Weil pairing.  Conversely, every $K$-defined maximal isotropic subgroup
$G$ of $(\Jac C)[2]$ arises in this way.  Thus, given a $K$-defined 
maximal isotropic subgroup $G$, we can define the {\em $G$-Richelot 
dual\/} of $C$ to be the Richelot dual of $C$ with respect to the 
factorization $f = g_1 g_2 g_3$ that gives rise to $G$.

\section{Proof of Theorem~\ref{T-nonsimple2}}
\label{S-nonsimple}

In this section we will prove Theorem~\ref{T-nonsimple2} by following 
the outline given in Section~\ref{S-curves} in the case where $A$ is a
split abelian surface.

Let $E$ and $E'$ be elliptic curves over a field $K$ of characteristic
not~$2$ and suppose there is a $2$-isogeny $\psi$ from $E$ to~$E'$.
Let $Q$ be the nonzero element of $E[2](K)$ in the kernel of $\psi$,
and let $P$ and $R$ be the other two geometric $2$-torsion points
of~$E$.  Let $Q' = \psi(P) = \psi(R)$, so that $Q'$ is a nonzero
element of $E'[2](K)$, and let $P'$ and $R'$ be the other geometric
$2$-torsion points of $E'$.  Suppose the discriminants of $E$ and $E'$
are equal up to squares, so that the fields $K(P)$ and $K(P')$ are the 
same.

Let $A$ be the surface $E\times E'$ and let $\lambda$ be the product 
principal polarization on~$A$.  Let $\alpha$ be the automorphism of 
$A$ that sends a point $(U,V)$ to $(U, V+\psi(U))$.  Let $G$ be the
$K$-defined subgroup $$G = \{(O,O), (P, P'), (Q, Q'), (R, R')\}$$
of $A[2]$ and let $G' = \alpha(G)$, so that $G'$ is the $K$-defined
subgroup $$G'= \{(O,O), (P, R'), (Q, Q'), (R, P')\}.$$  Let $(B,\mu)$ 
and $(B,\mu')$ be the principally-polarized surfaces obtained from $G$
and $G'$ as in Section~\ref{S-curves}.  The polarizations $\mu$ and
$\mu'$ will be indecomposable except in unusual circumstances, so 
there will usually be curves $C$ and $C'$ whose polarized Jacobians 
are isomorphic to $(B,\mu)$ and $(B,\mu')$, respectively.  If $E$ and 
$E'$ are given to us by explicit equations, then $C$ and $C'$ can also 
be given by explicit equations --- see~\cite[\S3.2]{HLP}, where the 
unusual circumstances are also explained.

To make this outline explicit and to prove Theorem~\ref{T-nonsimple2}
we must start with an explicit $2$-isogeny $\psi:E\to E'$, where the
discriminants of $E$ and $E'$ are equal up to squares.  Let $t$ be an
element of $K$ such that $t (t^2 + 1) (t^2 - 1)$ is nonzero, and let 
$E$ and $E'$ be the elliptic curves
\begin{align*}
E:  \quad y^2 & = x (x^2 - 4(t^2+1) x + 4(t^2+1))\\
E': \quad y^2 & = x (x^2 + 8(t^2+1) x + 16 t^2(t^2+1)).
\end{align*}
It is easy to check that the discriminants of $E$ and $E'$ are both 
equal to $t^2 + 1$, up to squares.

Let $s$ be a square root of $t^2 + 1$ in an algebraic closure of~$K$,
so that the $2$-torsion points of $E$ are
\begin{align*}
P &= (2t^2 + 2 + 2st, 0)\\
Q &= (0,0)\\
R &= (2t^2 + 2 - 2st,0)\\
\intertext{and the $2$-torsion points of $E'$ are}
P' &= (-4t^2 - 4 + 4s, 0)\\
Q' &= (0,0)\\
R' &= (-4t^2 - 4 - 4s,0).
\end{align*}
It is easy to check that the map
$$(x,y) \mapsto \left( \frac{y^2}{x^2}, 
                       \frac{(x^2 - 4(t^2+1))y}{x^2} \right)$$
defines a $2$-isogeny $\psi:E\to E'$ that kills $Q$ and that sends $P$ 
to~$Q'$ (see~\cite[Example III.4.5]{Silverman}).  Let $G$ and $G'$ be 
the subgroups of $A[2]$ defined above and let $(B,\mu)$ and $(B,\mu')$
be the principally-polarized surfaces obtained from $G$ and $G'$ as 
above.  If we apply~\cite[Prop.~4]{HLP} we find that $(B,\mu)$ is 
isomorphic over $K$ to the polarized Jacobian of the curve $y^2 = h_t$,
where 
$$h_t =  2^{38} t^6 (t+1)^3 (t^2+1)^{12}
         (2 x^2 - t)
         (4 t^2 x^4 + 4(t^2 + t + 1)x^2 + 1).$$
Furthermore, $(B,\mu')$ is isomorphic to the polarized Jacobians of the
curve $y^2 = h_{-t}$.  Scaling $h_t$ and $h_{-t}$ by squares in $K$, we
find that $y^2 = h_t$ is isomorphic to the curve $C(t)$ of
Theorem~\ref{T-nonsimple2} and that $y^2 = h_{-t}$ is isomorphic to the
curve~$C(-t)$.

To complete the proof of Theorem~\ref{T-nonsimple2} we must show that
$C(t)$ and $C(-t)$ are geometrically non-isomorphic, except for the
special cases listed in the theorem.  The simplest way to do this is to
use Igusa invariants (see~\cite{Igusa},~\cite{Mestre}).  Facilities for
computing Igusa invariants are included in the computer algebra
package Magma~\cite{Magma}.  

Let us begin by working over the ring $\BZ[t]$, where $t$ is an
indeterminate.   Let $J_2(t), J_4(t), J_6(t), J_8(t),$ and $J_{10}(t)$
be the Igusa invariants of the twist
$y^2 =  (2x^2 - t)(4x^4 + 4(t^2 + t + 1)x^2 + 1)$ of $C(t)$.  (The 
invariants $J_{2i}(t)$ of this curve, scaled by~$4^i$, can be computed
in Magma using the function {\tt ScaledIgusaInvariants}.)  Let
\begin{align*}
R_{2} &= \frac{J_4(t) J_2(-t)^2 - J_4(-t) J_2(t)^2}
              { t (t^2+1)^3 } \\
R_{3} &= \frac{J_6(t) J_2(-t)^3 - J_6(-t) J_2(t)^3}
              { t^3 (t^2+1)^3} \\
R_{5} &= \frac{J_{10}(t) J_2(-t)^5 - J_{10}(-t) J_2(t)^5}
              {t^3 (t^2+1)^7},
\end{align*}
all of which we view as elements of $\BZ[t]$.  If $C(t)$ and $C(-t)$ 
are isomorphic for a given value of $t$ in a given field~$K$, then the
polynomials $R_{2}, R_{3},$ and $R_{5}$ must all evaluate to $0$ at 
this value.  But we compute that 
$$\gcd(\Resultant(R_{2},R_{3}), \Resultant(R_{2},R_{5}))
  = 2^{980} 3^{48} 11^8,$$
so if the characteristic of $K$ is neither $3$ nor $11$ then the two 
curves $C(t)$ and $C(-t)$ are geometrically non-isomorphic for every 
value of $t$ in $K$ with $t(t^2+1)(t^2-1)\neq 0$.

We repeat the above calculation in the ring $\BF_3[t]$, only now we 
define
\begin{align*}
R_{2} &= \frac{J_4(t) J_2(-t)^2 - J_4(-t) J_2(t)^2}
              {t (t^2-1)^2 (t^2+1)^7} \\
R_{3} &= \frac{J_6(t) J_2(-t)^3 - J_6(-t) J_2(t)^3}
              {t^3 (t^2+1)^9}. 
\end{align*}
We find that $\gcd(R_{2},R_{3}) = 1$, so the two curves $C(t)$ and 
$C(-t)$ are geometrically non-isomorphic for every value of $t$ in 
characteristic~$3$, as long as $t(t^2+1)(t^2-1)\neq 0$.

Next we repeat the above calculation in the ring $\BF_{11}[t]$, with 
\begin{align*}
R_{2} &= \frac{J_4(t) J_2(-t)^2 - J_4(-t) J_2(t)^2}
              { t (t^2+1)^3 }\\
R_{3} &= \frac{J_6(t) J_2(-t)^3 - J_6(-t) J_2(t)^3}
              { t^3 (t^2+1)^3}. 
\end{align*}
We find that $\gcd(R_{2},R_{3}) = (t^2 + 3)(t^2 + 4)$, so that the two
curves are geometrically non-isomorphic in characteristic~$11$ except
possibly when $t^2$ is $-3$ or $-4$. 

Finally we note that in characteristic $11$, when $t^2$ is $-3$ or $-4$
the curve $C(t)$ is geometrically isomorphic to the supersingular curve
$y^2 = x^6 + x^4 + 4 x^2 + 7$.  Thus, $C(t)$ and $C(-t)$ are 
geometrically isomorphic for these values of $t$.
\qed

\begin{rem}
It was not critical in our construction that the isogeny $\psi:E\to E'$
have degree~$2$.  Similar constructions can be made with other kinds of
isogenies.
\end{rem}

\section{Jacobians with real multiplication by $\sqrt{2}$}
\label{S-RM}

In this section we review a construction of Bending~\cite{Bending} that
produces every genus-$2$ curve over a given field $K$ whose Jacobian
has a $K$-rational endomorphism that is fixed by the Rosati involution
and whose square is $2$.  We will give a variant of Bending's
construction that produces curves over $K$ with a not-necessarily
$K$-rational endomorphism that is fixed by Rosati and whose square
is~$2$.  We do not claim that our construction will produce all such
curves.

First we recall Bending's construction.  Let $K$ be a field of
characteristic not~$2$ and let $A$, $P$, and $Q$ be elements of $K$ 
with $P$ nonzero.  Define 
\begin{align*}
B &= (APQ - Q^2 + 4P^2 + 1)/P^2\\
C &= 4(AP - Q)/P\\
R &= 4P
\end{align*}
and let  $\alpha_1$, $\alpha_2$, $\alpha_3$ be the roots of
$T^3 + AT^2 + BT + C$ in a separable closure $\Kb$ of~$K$.  For 
$i = 1,2,3$ let 
$$G_i = X^2 - \alpha_i X + P\alpha_i^2 + Q\alpha_i + R,$$
and suppose that the product $G_1 G_2 G_3 \in K[X]$ has nonzero
discriminant. Let $D$ be a nonzero element of~$K$.

\begin{theorem}
\label{T-Bending}
The Jacobian of the genus-$2$ curve $D Y^2 = G_1 G_2 G_3$ has a
$K$-rational endomorphism that is fixed by the Rosati involution and 
whose square is $2$.  Furthermore, if $\#K>5$ then every curve over $K$
whose Jacobian has such an endomorphism is isomorphic to a curve that
arises in this way from some choice of $A$, $P$, $Q$, and $D$ in~$K$.
\end{theorem}

\begin{proof}
See~\cite[Theorem 4.1]{Bending}.  Bending assumes that the base field
$K$ has characteristic~$0$, but his proof works over an arbitrary field
$K$ of characteristic not~$2$ so long as every genus-$2$ curve over $K$
can be written in the form $y^2 = \text{(sextic)}$.  This is the case 
for every field with more than $5$ elements.
\end{proof}

The endomorphism of $D Y^2 = G_1 G_2 G_3$ whose existence is claimed
by Theorem~\ref{T-Bending} is obtained by noting that the obvious 
Richelot dual of the curve is isomorphic over $K$ to the curve itself.
Thus the degree-$4$ isogeny from the Jacobian of the curve to the
Jacobian of its dual can be viewed as an endomorphism of the curve's
Jacobian, and this endomorphism has the properties claimed in the
theorem.

We will want to consider curves over $K$ whose Jacobians have real 
multiplication by $\sqrt{2}$ that is not necessarily defined over~$K$.
For this reason, we will require the following variant of Bending's 
construction:

Suppose $r$, $s$, and $t$ are elements of a field $K$ of characteristic
not~$2$, with $s\neq 0$, $s\neq 1$, and $t\neq 1$.  Let 
\begin{align*}
c_2 &= r + 4t\\
c_1 &= 4 t (r + s^3 - s^2 t - 2 s^2 + 5 s + t)\\
c_0 &= 4 t (s - 1) (r s^2 - r s t - r s - r t - 8 s t)
\end{align*}
and suppose that the polynomial $T^3 - c_2 T^2 + c_1 T - c_0$ has three
distinct roots $\beta_1$, $\beta_2$, $\beta_3$ in $\Kb$.  For 
$i = 1,2,3$ let 
$$g_i = x^2 - 
        2\beta_i x + 
        (1-s)\beta_i^2 - 4 s (s - 1)^2 t (s - t - 1),$$
and suppose that the discriminant of the product $f = g_1 g_2 g_3$ is
nonzero. Let $\CC(r,s,t)$ be the curve over $K$ defined by $y^2 = f$.

\begin{theorem}
\label{T-RM}
The Richelot dual of $\CC(r,s,t)$ associated to the factorization 
$f = g_1 g_2 g_3$ is isomorphic over $K(\sqrt{st}\,)$ to~$\CC(r,s,t)$. 
The endomorphism of $\Jac \CC(r,s,t)$ {\rm(}over 
$K(\sqrt{st}\,)$\/{\rm)} obtained by composing the Richelot isogeny
with the natural isomorphism from the dual curve to $\CC(r,s,t)$ is 
fixed by the Rosati involution, and its square is the
multiplication-by-$2$ endomorphism.
\end{theorem}

\begin{proof}
This theorem can be proven by direct calculation, but here we will
prove it by relating the curve $\CC(r,s,t)$ back to Bending's 
construction.  Let $Q$ be a square root of $st$ in $\Kb$, let 
$P = (1 - s)/4$, let $A = (r + 6 s t - 2 t)/(4PQ)$, let $D=1$, and let
$\CC'$ be the curve over $K(Q)$ defined by using this $A$, $P$, $Q$, 
and $D$ in Bending's construction.  Then the $\alpha_i$ are related to
the $\beta_i$ by $$Q\alpha_i = \beta_i/(s-1) + 2t,$$ and the curve 
$Y^2 = G_1 G_2 G_3$ is isomorphic $y^2 = g_1 g_2 g_3$ via the relation
$x = 2(s-1)(QX + t)$.  This shows that $\CC(r,s,t)$ is isomorphic to
its Richelot dual over $K(Q)$.  The rest of the theorem follows from 
Bending's theorem.
\end{proof}

\begin{rem}
Bending's family of curves has three ``geometric'' parameters $A$, $P$,
and $Q$ and one ``arithmetic'' parameter $D$ (which parametrizes
quadratic twists of the curve determined by $A$, $P$, and $Q$).  Since
the moduli space $\CM$ of genus-$2$ curves with real multiplication by
$\sqrt{2}$ is a two-dimensional rational variety, one might hope to
replace Bending's three-geometric-parameter family with a two-parameter
family.  But there is an obstruction, which stems from the fact that
$\CM$ is a coarse moduli space and not a fine one:  A $K$-rational
point on $\CM$ does not necessarily give rise to a curve over $K$.
Indeed, Mestre~\cite{Mestre} has shown that to every $K$-rational point
$P$ on the moduli space of genus-$2$ curves there is naturally 
associated a genus-$0$ curve over $K$, and $P$ corresponds to a curve
over $K$ if and only if the genus-$0$ curve has a $K$-rational point.  
\end{rem}

\section{Galois restrictions}
\label{S-Galois}

In order to prove Theorem~\ref{T-simple2} we will apply the
construction outlined in Section~\ref{S-curves} to a Jacobian with real
multiplication by~$\sqrt{2}$ that we will obtain from 
Theorem~\ref{T-RM};  we will take the automorphism $\alpha$ to be 
$1+\sqrt{2}$.  The construction requires that we find a Galois-stable 
maximal isotropic subgroup $G$ of the $2$-torsion of the Jacobian such 
that $G' = (1+\sqrt{2})(G)$ is a maximal isotropic subgroup different 
from~$G$.  This requirement imposes some restrictions on the values of
$r$, $s$, and $t$ that we will be able to use in Theorem~\ref{T-RM}.  
In this section we will make these restrictions explicit, and in 
Section~\ref{S-application} we will find an elliptic surface that
parametrizes a subset of the allowable values of $r$, $s$, and~$t$.

Recall the basic outline of Theorem~\ref{T-RM}: Given three elements 
$r$, $s$, $t$ of our base field $K$, we define a polynomial
$h = T^3 - c_2 T^2 + c_1 T - c_0$ in the polynomial ring $K[T]$, and we
assume that $h$ is separable.  We use the roots $\beta_1$, $\beta_2$,
$\beta_3$ of $h$ to define three polynomials $g_1$, $g_2$, $g_3$ in the
polynomial ring $\Kb[x]$, we assume that the product $f = g_1 g_2 g_3$
is separable, and we define a curve $C$ by $y^2 = f$.  Then we show
that the Richelot dual of $C$ corresponding to the factorization
$f = g_1 g_2 g_3$ is geometrically isomorphic to $C$ itself.

Let $L$ be the quotient of the polynomial ring $K[T]$ by the ideal
generated by the polynomial $h$ and let $\beta$ be the image of $T$ 
in $L$.  Since $h$ is separable, the algebra $L$ is a product of 
fields.  Let $g\in L[x]$ be the polynomial
$$g = x^2 - 2\beta x + (1-s)\beta^2 - 4 s (s-1)^2 t (s - t - 1).$$
Let $\Delta\in K^*$ be the discriminant of $h$ and let $\Delta'\in L^*$
be the discriminant of~$g$.

\begin{theorem}
\label{T-Galois}
There are distinct Galois-stable maximal isotropic subgroups $G$ and 
$G'$ of $(\Jac C)[2]$ with $G' = (1+\sqrt{2})(G)$ if and only if 
$\Delta\Delta'$ is a square in the algebra~$L$.
\end{theorem}

\begin{proof}
Let the roots of $g_1$ (respectively, $g_2$, $g_3$) be $r_1$ and $r_2$
(respectively, $r_3$ and $r_4$, $r_5$ and $r_6$).  For each $i$ let
$W_i$ be the Weierstra\ss\ point of $C$ corresponding to the root $r_i$
of $f = g_1 g_2 g_3$.  The kernel $H$ of the Richelot isogeny
multiplication-by-$\sqrt{2}$ on the Jacobian $J$ of $C$ is the 
order-$4$ subgroup containing the divisor classes $[W_1 - W_2]$, 
$[W_3 - W_4]$, and $[W_5 - W_6]$.

Suppose there are distinct Galois-stable maximal isotropic subgroups 
$G$ and $G'$ of $(\Jac C)[2]$ with $G' = (1+\sqrt{2})(G)$.  Then 
clearly $G\neq H$, so $\#(G\cap H)$ is either $2$ or~$1$.  Suppose 
$\#(G\cap H)=2$.  By renumbering the polynomials $g_i$ and by 
renumbering their roots, we may assume that $G$ is the order-$4$
subgroup $$G = \{ 0, [W_1 - W_2], [W_3 - W_5], [W_4 - W_6]\}.$$
Then $\sqrt{2}$ kills the first two elements of $G$ and sends the 
second two elements to $[W_1 - W_2]$, and it follows that 
$(1+\sqrt{2})(G) = G$, contradicting our assumption that $G$ and $G'$
are distinct.

So now we know that $G\cap H = \{0\}$.  By renumbering the polynomials
$g_i$ and by renumbering their roots, we may assume that $G$ is the 
order-$4$ subgroup 
$$G = \{ 0, [W_1 - W_5], [W_2 - W_4], [W_3 - W_6]\}.$$  It is not hard
to show that the automorphism $1 + \sqrt{2}$ of $J$ sends $[W_1 - W_5]$
to $[W_2 - W_6]$, $[W_2 - W_4]$ to $[W_1 - W_3]$, and $[W_3 - W_6]$ to
$[W_4 - W_5]$, so we have 
$$G'= \{ 0, [W_1 - W_3], [W_2 - W_6], [W_4 - W_5]\}.$$

Suppose $\sigma$ is an element of the Galois group such that 
$r_1^\sigma = r_2$.  Since $G$ is Galois stable, it follows that
$\sigma$ sends $[W_1 - W_5]$ to $[W_2 - W_4]$, and therefore 
$r_5^\sigma = r_4$.  But since $G'$ is Galois stable, we see that 
$\sigma$ must send $[W_4 - W_5]$ to itself, and it follows that 
$r_4^\sigma = r_5$.  Continuing in this manner, we find that 
$r_2^\sigma = r_1$ and $r_6^\sigma = r_3$ and $r_3^\sigma = r_6$.
Thus, $\sigma$ acts on the roots of $f$ according to the permutation 
$(1 2)(3 6)(4 5)$ of the subscripts.

By considering the other choices for $r_1^\sigma$ and using the same 
reasoning as above, we find that the image of the absolute Galois group
of $K$ in the symmetric group on the the roots of $f$ is contained in
the subgroup
$$S = \{ \Id, (1 2)(3 6)(4 5), (1 3)(2 4)(5 6), (1 4 6)(2 3 5),
   (1 5) (2 6) (3 4), (1 6 4) (2 5 3)\};$$
here of course we identify the root $r_i$ with the integer~$i$.  In 
particular, note that the action of $\sigma$ on the $r_i$ is determined
by the action of $\sigma$ on the~$\beta_i$.

To show that $\Delta\Delta'$ is square in the algebra $L$, we will
consider three cases, depending on the splitting of the polynomial~$h$.

\emph{Case 1.}  Suppose $h$ is irreducible.  Then $L$ is a field, and
the condition that $\Delta\Delta'$ be a square in $L$ is equivalent to 
saying that $g$ defines the Galois closure $M$ of $L$ over~$K$.  So 
suppose, to obtain a contradiction, that $g$ does not define $M$ 
over~$L$.  There are two ways that this can happen: either the roots 
of $g$ do not lie in $M$, or $M$ is a quadratic extension of $L$ and 
the roots of $g$ lie in~$L$.

Suppose that the roots of $g$ do not lie in~$M$.  Then there is an 
element $\sigma$ of the absolute Galois group of $K$ that fixes $M$ 
but that moves the roots of~$g$.  But this contradicts the fact that 
the the action of $\sigma$ on the roots of $g$ is determined by the 
action of $\sigma$ on the~$\beta_i$.

Suppose $M$ is a quadratic extension of $L$, so that  the image of the
absolute Galois group in the symmetric group on the $\beta_i$ is the
full symmetric group.  Then the image of the absolute Galois group in
the symmetric group on the $r_i$ must be the entire group $S$ given 
above, which acts transitively on the~$r_i$.  But if the roots of $g$
lie in $L$, then the $r_i$ will form two orbits under the action of the
absolute Galois group, giving a contradiction.

\emph{Case 2.}  Suppose $h$ factors as a linear polynomial times an 
irreducible quadratic.  Then one of the $\beta_i$, say $\beta_1$, lies
in $K$, while $\beta_2$ and $\beta_3$ are conjugate elements in a
quadratic extension $M$ of~$K$.  With the labelings we have chosen,
this means that the image of the absolute Galois group in the symmetric
group on the $r_i$ must be equal to the two-element group
$$S' = \{ \Id, (1 2)(3 6)(4 5)\};$$  In particular, we see that $r_3$
and $r_6$ (and $r_4$ and $r_5$) are quadratic conjugates of one 
another, so they all must be elements of~$M$.  This means that the 
image of $g$ in $K[x]$ (obtained by sending $\beta$ to~$\beta_1$) is 
an irreducible polynomial that defines $M$, while the image of $g$ in
$M[x]$ (obtained by sending $\beta$ to~$\beta_2$) splits into two 
linear factors.  Thus, the discriminant $\Delta'$ of $g$ in
$L = K\times M$ is equal to the discriminant of $M$ (up to squares) in
the first component, and is a square in the second component.  But 
$\Delta$ has this same property, so the product $\Delta\Delta'$ is a
square in~$L$.

\emph{Case 3.}  Suppose $h$ splits over $K$ into three linear factors,
so that $\Delta$ is a square in~$K$.  Then the absolute Galois group of
$K$ acts trivially on the $\beta_i$, so it must act trivially on the
$r_i$ as well.  This means that the discriminant $\Delta'$ of $g$ must
be a square in each factor of $L = K\times K\times K$, so 
$\Delta\Delta'$ is a square in $L$ as well.

We see that if there are subgroups $G$ and $G'$ as in the statement of
the theorem then $\Delta\Delta'$ must be a square in~$L$.  We leave the
details of the proof of the converse statement to the reader; the point
is that in each of the three cases above, the reasoning is reversible.
\end{proof}

\section{Application of our construction to \\
curves with real multiplication}
\label{S-application}

In this section we will follow the outline given in 
Section~\ref{S-curves} in the case where $A$ is a Jacobian with real
multiplication by $\sqrt{2}$ that has appropriate Galois-stable 
subgroups.  Theorem~\ref{T-simple2} will follow quickly from the result
we obtain.

We will continue to use the notation from previous sections:
\begin{itemize}
\item
$r$, $s$, and $t$ will be elements of a field $K$;
\item
$h$ will be a polynomial in $K[T]$ defined in terms of 
$r$, $s$, and $t$; 
\item
$L$ will be the algebra $K[T]/(h)$;
\item
$\beta$ will be the image of $T$ in $L$;
\item
$g$ will be a polynomial in $L[x]$ defined in terms
of $r$, $s$, $t$, and $\beta$;
\item
$\Delta\in K^*$ will be the discriminant of $h$; and
\item
$\Delta' \in L^*$ will be the discriminant of $g$.
\end{itemize}

\begin{theorem}
\label{T-gensimple}
Let $K$ be a field of characteristic not~$2$, suppose $r$, $s$, and $t$
are elements of~$K$ that satisfy the hypotheses appearing before the 
statement of Theorem~{\rm\ref{T-RM}}, and let $C$ be the curve 
$\CC(r,s,t)$ from  Theorem~{\rm\ref{T-RM}}.  Suppose further that the
product $\Delta\Delta'$ is a square in~$L^*$, so that there are 
Galois-stable subgroups $G$ and $G'$ of $(\Jac C)[2]$ as in the 
statement of Theorem~{\rm\ref{T-Galois}}.  Then the Jacobian of the 
$G$-Richelot dual of $C$ is isomorphic over $K(\sqrt{st}\,)$ to the 
Jacobian of the $G'$-Richelot dual of~$C$.
\end{theorem}

\begin{proof}
Let $D$ be the $G$-Richelot dual of $C$ and let $D'$ be the 
$G'$-Richelot dual of~$C$.  Theorem~\ref{T-RM} and 
Theorem~\ref{T-Galois}, combined with the argument in 
Section~\ref{S-curves}, show that the Jacobian of $D$ becomes 
isomorphic to the Jacobian of $D'$ when the base field is extended 
to~$K(\sqrt{st}\,)$.
\end{proof}

\begin{rem}
Since $D$ and $D'$ are defined over $K$, and since their Jacobians
become isomorphic over $K(\sqrt{st}\,)$, it is tempting to think that
$\Jac D$ must be isomorphic over $K$ to either the Jacobian of $D'$ or
the Jacobian of the standard quadratic twist of $D'$ 
over~$K(\sqrt{st}\,)$.  But in fact this is not the case.  It is true
that $\Jac D'$ is a $K(\sqrt{st}\,)/K$-twist of $\Jac D$, but the twist
is by an automorphism of $\Jac D$ that does not come from an
automorphism of~$D$.  Indeed, generically the automorphism group of $D$
contains $2$ elements, while the automorphism group of $\Jac D$ is 
isomorphic to the unit group of $\BZ[2\sqrt{2}]$.
\end{rem}

\begin{proposition}
Suppose $t = s-1$ and let $u = r + 2$.  Then $\Delta\Delta'$ is a 
square in $L$ if and only if $(u^2 + a)^2 + 8bu + 4c$ is a square 
in~$K$, where
\begin{align*}
a &= -4 s (s^2 + 11s - 11)\\
b &= -8s^2 (s-1) (4s-1) \\
c &= -16s^2 (s-1) (28s^2 - 19s + 1).
\end{align*}
\end{proposition}

\begin{proof}
When $t = s-1$ and $r = u - 2$, we find that the coefficients of the 
polynomial $h$ used to define the algebra $L$ are
\begin{align*}
c_2 &= 4s + u - 6 \\
c_1 &= - 4 (s-1) (s^2 - 6s - u + 3) \\
c_0 &= 4 (s-1)^3 (-8s - u + 2),
\end{align*}
and we compute that $$\Delta = 16 s (s-1)^2 ((u^2 + a)^2 + 8bu + 4c),$$
where $a$, $b$, and $c$ are as in the statement of the proposition.
Furthermore, the polynomial $g\in L[x]$ defined in 
Section~\ref{S-Galois} is $x^2 - 2\beta x + (1-s)\beta^2$, so that 
$\Delta' = 4 s \beta^2$.  We see that $\Delta\Delta'$ is a square in 
$L$ if and only if the element $\delta = (u^2 + a)^2 + 8bu + 4c$ of $K$
is a square in~$L$.  If $L$ is a field then it is a cubic extension
of~$K$, and $\delta$ is a square in $L$ if and only if it is a square 
in~$K$.  If $L$ is not a field then it has $K$ as a factor, and again
$\delta$ is a square in $L$ if and only if it is a square in~$K$.
\end{proof}

\begin{proposition}
\label{P-ellipticsurface}
Let $K=\BQ(s)$ be the function field in the variable $s$ over~$\BQ$, 
let 
\begin{align*}
a &= -4 s (s^2 + 11s - 11)\\
b &= -8s^2 (s-1) (4s-1) \\
c &= -16s^2 (s-1) (28s^2 - 19s + 1),
\end{align*}
let $F$ be the curve over $K$ defined by 
$$z^2 = (u^2 + a)^2 + 8bu + 4c,$$ and let $E$ be the elliptic curve
over $K$ defined by $$y^2 = x^3 - a x^2 - c x + b^2.$$  Then 
\begin{itemize}
\item[(a)] the map $u = (y-b)/x$, $z = 2x - u^2 - a$ gives 
           an isomorphism from $E$ to $F$, whose
           inverse is $x = (z + u^2 + a)/2$, $y = ux + b${\/\rm;}
\item[(b)] the point $P = (0, b)$ on $E$ has infinite order{\/\rm;}
\item[(c)] the point $T = (4s^2(1-s), 0)$ on $E$ has order $2${\/\rm;}
\item[(d)] the isomorphism in statement {\rm(a)}
           takes the involution $(u,z)\mapsto (u,-z)$ on $F$
           to the involution $Q \mapsto -Q-P$ on $E${\/\rm;}
\item[(e)] the isomorphism in statement {\rm(a)} takes $-P$
           and the origin of $E$ to the two infinite points on $F$.
\end{itemize}
\end{proposition}             

\begin{proof}
An easy calculation shows that statement~(a) is true; the particular 
values of $a$, $b$, and $c$ are irrelevant to the calculation.

To show that the point $P$ has infinite order, it suffices to show that
when we specialize $s$ to a particular value the specialized $P$ has
infinite order.  For example, if we set $s = 2$, then $E$ becomes the
curve $$y^2 = x^3 + 120 x^2 + 4800 x + 50176$$ and $P$ becomes the
point $(0, -224)$.  Translating $x$ by $40$, we find a new equation 
for $E$: $y^2 = x^3 - 13824,$ where now $P = (40,-224)$.  But $7$
divides $224$ and $7$ does not divide $13824$, so by the Lutz-Nagell
theorem~\cite[Cor.~VIII.7.2]{Silverman} the point $P$ has infinite 
order.  This proves statement~(b).

Statement~(c) is clear.

Let $R$ be a point $(u,z)$ on $F$ and let $\tilde{R} = (u,-z)$ be its
involute.  Let $Q$ and $\tilde{Q}$ be the images of $R$ and $\tilde{R}$
on $E$.  Clearly $Q$ and $\tilde{Q}$ both lie on the line $y = ux + b$,
and the third intersection point of this line with $E$ is easily seen 
to be~$P$.  Thus, the involution on $E$ satisfies $\tilde{Q} + Q = -P$,
and this is statement~(d).

The equations for the isomorphism show that $-P$ is mapped to an 
infinite point on $F$, and statement~(d) shows that $O_E$ gets mapped
to an infinite point as well.
\end{proof}

\begin{rem}
If we view the curve $F\cong E$ as an elliptic surface $S$ over~$\BQ$, 
then the points $P$ and $T$ of Proposition~\ref{P-ellipticsurface} can
be viewed as rational curves on~$S$.  By adding multiples of $P$ and 
$T$ together, we get a countable family of rational curves on~$S$.  But
$S$ contains more rational curves than just the ones in this family.
For example, we have the curves
\begin{align*}
s &= 5/4 \\
u &= (4w^2 + 5w + 40)/(4w) \\
z &= (2w^4 + 5w^3 - 50w - 200)/(2w^2)
\end{align*}
and
\begin{align*}
s &= -1 \\
u &= (2w^2 - 10w - 4)/w \\
z &= (4w^4 - 40w^3 - 80w - 16)/w^2,
\end{align*}
where $w$ is a parameter; the curve
\begin{align*}
s &= (5-w^2)/4\\
u &= (-w^4 + w^3 + 7w^2 - 5w - 10)/(4w + 8) \\
z &= (w^8 + 9w^7 + 22w^6 - 18w^5 - 135w^4 - 135w^3)/(8w^2 + 32w + 32),
\end{align*}
which corresponds to a $3$-torsion point on $E$ defined over a 
genus-$0$ extension of the function field~$\BQ(s)$;  five curves in
which $u$ is a linear expression in~$s$, for example
\begin{align*}
s &= 4(w^2 + 9w + 19)/w\\
u &= -s\\
z &= 16(16w^6 + 283w^5 + 1555w^4 - 29545w^2 - 102163w - 109744)/w^3;
\end{align*}
and three curves in which $u$ is a quadratic expression in $s$, for
example
\begin{align*}
s &= (w^2 + 3w + 1)/w\\
u &= 4s^2 - 6s\\
z &= 4(4w^8 + 35w^7 + 105w^6 + 119w^5 - 119w^3 - 105w^2 - 35w - 4)/w^4.
\end{align*}
\end{rem}

\begin{rem}
One can check that the image of the elliptic surface $S$ in the moduli
space of genus-$2$ curves is $2$-dimensional.  To check this, one need
only write explicitly the Igusa invariants of the genus-$2$ curve
obtained from a pair $(s,u)$ and verify that the rank of the Jacobian
matrix of the mapping from $(s,u)$-pairs to Igusa invariants at some
arbitrary point is~$2$.
\end{rem}

We now have enough machinery available to prove 
Theorem~\ref{T-simple2}.

\begin{proof}[Proof of Theorem~{\rm\ref{T-simple2}}]
Consider the point $P = (0,b)$ from statement~(b) of
Proposition~\ref{P-ellipticsurface}.  The $u$-co\"ordinate of its image
on the curve $F$ is $$u = (28s^2 - 19s + 1)/(1-4s).$$  So given any 
$s\in K$, we will obtain a curve satisfying the conclusion of 
Theorem~\ref{T-gensimple} if we set
\begin{align*}
 t &= s - 1\\
 r &= -2 + (28s^2 - 19s + 1)/(1-4s) 
\end{align*}
and set $C = \CC(r,s,t)$.  Given a $v \in K\setminus\{0,1,4\}$, let us
apply the preceding observation with $s = v/4$, so that 
\begin{align*}
s &= v/4\\
t &= (v - 4)/4\\
r &= (7 v^2 - 11 v - 4)/(4-4v).
\end{align*}
The coefficients of the polynomial $h$ used in the construction of 
Section~\ref{S-RM} are
\begin{align*}
c_2 &= \frac{3v^2 + 9v - 20}{4(1 - v)}             \\
c_1 &= \frac{(v-4)(v^3 + 3v^2 - 4v - 32)}{16(1-v)}  \\
c_0 &= \frac{(v-4)^3 (v^2 + 3v + 4)}{64(1-v)},
\end{align*}
and over $K(w)$ the roots of $h$ are
\begin{align*}
\beta_1 &= \frac{(2 + w)(2 - w)}{4} \\
\beta_2 &= \frac{-(2+w)^2(2 - w + w^2)}{4(1 + w)}\\
\beta_3 &= \frac{-(2-w)^2(2 + w + w^2)}{4(1 - w)}.
\end{align*}
Each polynomial $g_i$ is $x^2 - 2\beta_i x + (1-s)\beta_i^2$ and has 
roots $\beta_i (1 \pm w/2)$, so we calculate that the roots of $g_1$
are
\begin{align*}
r_1 &= -(1/8) (2 - w)^2 (2 + w)                           \\
r_2 &= -(1/8) (2 - w)   (2 + w)^2,                        \\
\intertext{the roots of $g_2$ are}
r_3 &= -(1/8)           (2 + w)^3 (2 - w + w^2) / (1 + w) \\
r_4 &= -(1/8) (2 - w)   (2 + w)^2 (2 - w + w^2) / (1 + w),\\
\intertext{and the roots of $g_3$ are}
r_5 &= -(1/8) (2 - w)^2 (2 + w)   (2 + w + w^2) / (1 - w) \\
r_6 &= -(1/8) (2 - w)^3           (2 + w + w^2) / (1 - w).
\end{align*}
These roots are indexed in a manner consistent with the indexing of the
roots in the proof of Theorem~\ref{T-Galois}.  Note that the $r_i$ are
related to the $\rho_i$ of Theorem~\ref{T-simple2} by the relation
$$r_i = 4s(s-1)\rho_i - 2(s-1)^2,$$
so the $r_i$ are distinct exactly when the $\rho_i$ are distinct.  It
is easy to check that when 
$(v^2 + 3 v + 4)(v^2 - v + 4) (v^3 - 6 v^2 - 7 v - 4)\neq 0$ the 
$\rho_i$ are distinct, so in this case the curve $D$ of 
Theorem~\ref{T-simple2} has genus~$2$.  Furthermore, we see that $D$ is
isomorphic to the curve $\CC(r,s,t).$

For each $i$ let $W_i$ be the point $(\rho_i,0)$ of $D$.  Let $G$ be
the Galois-stable subgroup of the Jacobian of $D$ that consists of the
divisor classes $$\{ [0], [W_1 - W_5], [W_2 - W_4], [W_3 - W_6]\}$$
and let $G'$ be the Galois-stable subgroup 
$$\{ [0], [W_1 - W_3], [W_2 - W_6], [W_4 - W_5]\}.$$ 
An easy computation shows that when $v^3 - 4 v^2 + 7 v + 4$ and 
$v^2 + v + 2$ are nonzero the $G$-Richelot dual of $D$ and the 
$G'$-Richelot dual of $D$ are defined (that is, the determinants
mentioned in the definition of the two Richelot duals are nonzero).
Then the results of Section~\ref{S-RM} show that the $G$-Richelot dual
of $D$ and the $G'$-Richelot dual of $D$ become isomorphic over 
$K(\sqrt{st}\,) = K(\sqrt{v(v-4)}\,)$.

The proof that these two Richelot duals of $D$ are geometrically 
non-isomorphic to one another (except in the special cases listed in 
the theorem) is a computation along the same lines as the proof of the
corresponding statement of Theorem~\ref{T-nonsimple2}.  We leave the
details to the reader.

Finally, suppose that $K$ has characteristic $0$ and suppose that
there is a ring homomorphism from $\BZ[v]$ to $\BF_{13}$ that takes
$v$ to either $2$ or~$6$.  Then the curve $D$ reduces modulo $13$ to
one of the curves over $\BF_{13}$ obtained when $v=2$ or~$v=6$.  One 
can compute that the characteristic polynomials of Frobenius for the 
Jacobians of these two curves are $t^4 - 4 t^3 + 22 t^2 - 52 t + 169$
and $t^4 + 4 t^3 + 22 t^2 + 52 t + 169$, respectively.  
Then~\cite[Thm.~6]{Howe-Zhu} shows that the Jacobians are absolutely 
simple.  Since $D$ modulo $13$ is absolutely simple, so is $D$ itself.
And finally, since $C$ and $C'$ have Jacobians isogenous to that
of~$D$, we see that their Jacobians are absolutely simple too.  The
final statement of the theorem then follows from the observation that
if $v$ is not algebraic, then there is a homomorphism 
$\BZ[v]\to\BF_{13}$ that sends $v$ to any given element.
\end{proof}

\begin{rem}
We used the field $\BF_{13}$ at the end of the proof simply because it 
is the smallest prime field that contains values of $v$ that give rise 
to absolutely simple Jacobians.  Other prime fields have a larger 
proportion of good values of $v$.  For example, there are $341$ values
of $v$ in $\BF_{769}$ that give rise to absolutely simple Jacobians.
For three-digit primes $p$ the number of good $v$-values is typically 
greater than~$0.3 p$.  This implies that for a ``randomly chosen''
rational number $v$, it is almost certainly the case that $v$ will give
rise to an absolutely simple Jacobian.
\end{rem}

\section{Proof of Theorem~{\rm\ref{T-nonsimple3}}}
\label{S-genus3}

The proof of Theorem~\ref{T-nonsimple3} is very much like the proof of
Theorem~\ref{T-nonsimple2}:  We will produce three elliptic curves 
$E_1$, $E_2$, $E_3$, two maximal isotropic subgroups $G$, $G'$ of the 
$2$-torsion of $A = E_1\times E_2\times E_3$, and an automorphism
$\alpha$ of $A$ that takes $G$ to~$G'$.  Then we will produce a 
hyperelliptic curve whose Jacobian is $A/G$ and a plane quartic whose
Jacobian is $A/G'$.  To produce these curves we will use the results 
of~\cite[\S4]{HLP}.  Our notation will be chosen to match that 
of~\cite{HLP}, except that we will continue to call our base field~$K$,
instead of~$k$.

Let $K$ be an arbitrary field of characteristic not $2$ and let $t$ be
an element of $K$ with $t (t + 1) (t^2 + 1) (t^2 + t + 1) \neq 0$.  Let
$s = -(t^2 + t + 1)$ and let $r$ be a square root of $t^2 + 1$ in an 
algebraic closure of~$K$.  Let
\begin{center}
\renewcommand{\arraystretch}{1.5}
\begin{tabular}{lll}
$A_1 = -2(t^2+1)s$         & $B_1 =     (t^2+1)s^2$   \\
$A_2 =  4(t^2+1)s$         & $B_2 = 4t^2(t^2+1)s^2$   \\
$A_3 = -2(t^2+t+1)s\qquad$ & $B_3 = (t+1)^2(t^2+1)s^2$\\
\end{tabular}
\end{center}
and for each $i$ let 
\begin{align*}
\Delta_1 &= A_1^2 - 4B_1 = 4t^2(t^2+1)s^2 \\
\Delta_2 &= A_2^2 - 4B_2 = 16(t^2+1)s^2 \\
\Delta_3 &= A_3^2 - 4B_3 = 4t^2s^2.
\end{align*}
Note that the $\Delta_i$ and the $B_i$ are nonzero, so we may define
for each $i$ an elliptic curve $E_i$ by $$y^2 = x(x^2 + A_i x + B_i).$$
We define $2$-torsion points $P_i$ on the $E_i$ by setting
\begin{align*}
P_1 &= \bigl((t^2+1)s - rts,   0\bigr)\\
P_2 &= \bigl( -2(t^2+1)s - 2rs, 0\bigr)\\
P_3 &= \bigl((t^2+t+1)s - ts,  0\bigr)
\end{align*}
and for each $i$ we let $Q_i$ be the $2$-torsion point $(0,0)$ on~$E_i$
and we let $R_i = P_i + Q_i$.

Let $A = E_1\times E_2\times E_3$ and let $G$ be the subgroup of $A[2]$
generated by $(P_1,P_2,P_3)$, $(Q_1,Q_2,0)$, and $(Q_1,0,Q_3)$.  
Associated to these choices of $A$ and $G$ there is a quantity called
the {\it twisting factor\/} $T$ (see~\cite[\S4]{HLP}).  Using the
formula in~\cite[\S4]{HLP} we find that for our $A$ and $G$ the 
twisting factor is~$0$, so we may apply~\cite[Prop.~14]{HLP} to find a 
hyperelliptic genus-$3$ curve whose Jacobian is isomorphic over $K$ 
to~$A/G$.  The curve given by~\cite[Prop.~14]{HLP} is defined by two
equations in~$\BP^3$, namely 
\begin{align}
W^2 Z^2 &= a X^4 + b Y^4 + c Z^4 \label{E-oldH1}\\
      0 &= d X^2 + e Y^2 + f Z^2 \label{E-oldH2}
\end{align}
where
\begin{align*}
  a &=  4 t   (t+1)   (t^2+1)^3 s^5 \\
  b &= 16 t^2 (t+1)   (t^2+1)^3 s^5 \\
  c &=  4 t   (t+1)^2 (t^2+1)^2 s^6 \\
1/d &= -2 t   (t+1)   (t^2+1) s^2   \\
1/e &=        (t+1)   (t^2+1) s^2   \\
1/f &=  2 t           (t^2+1) s^2.
\end{align*}
If we replace $W$ by $2t(t+1)(t^2+1)s^3 W$ in Equation~\ref{E-oldH1}
and divide out common factors, we get Equation~\ref{E-H1}, and if we
multiply Equation~\ref{E-oldH2} by $2 t (t+1) (t^2+1) s^2$ we get 
Equation~\ref{E-H2}. This shows that the Jacobian of $H(t)$ is 
isomorphic to $A/G$.

Now let $G'$ be the subgroup of $A[2]$ generated by $(P_1,P_2,R_3)$, 
$(Q_1,Q_2,0)$, and $(Q_1,0,Q_3)$, and let $T'$ be the twisting factor
associated to $A$ and $G'$.  The formula in~\cite[\S4]{HLP} shows that
$$T' = -64(t^2+1)^2(t^2+t+1)s^3 = 64(t^2+1)^2s^4,$$ so the twisting
factor is a nonzero square.  Then~\cite[Prop.~15]{HLP} shows that there
is a plane quartic whose Jacobian is isomorphic (over~$K$) to~$A/G'$.  
The plane quartic is given by
\begin{equation}
\label{E-oldQ}
B_1 X^4 + B_2 Y^4 + B_3 Z^4 + d' X^2Y^2 + e' X^2Z^2 + f'Y^2Z^2 = 0
\end{equation}
where
\begin{align*}
d' &=  4 (t^2 + 1) (2t^2 + t + 2) s^2\\
e' &= -2 (t^2 + 1) (2t^2 + t + 1) s^2\\
f' &=  4 (t^2 + 1) ( t^2 + t + 2) s^2.
\end{align*}
Dividing Equation~\ref{E-oldQ} by $(t^2+1)s^2$ gives 
Equation~\ref{E-Q}, so the Jacobian of $Q(t)$ is isomorphic to $A/G'$.

To complete the proof we must show that $A/G \cong A/G'$.  Note that 
there is a $2$-isogeny $\psi$ from $E_1$ to $E_2$ that kills $Q_1$ and
that takes $P_1$ and $R_1$ to $Q_2$ 
(see~\cite[Example III.4.5]{Silverman}).  Consider the automorphism 
$\alpha$ of $A$ that sends a point $(S_1,S_2,S_3)$ to 
$(S_1,S_2+\psi(S_1),S_3)$.  It is easy to check that $\alpha(G) = G'$,
and it follows that $A/G\cong A/G'$, as desired, so the Jacobians of
$H(t)$ and $Q(t)$ are isomorphic over~$K$.
\qed

\section{Examples}
\label{S-examples}

\begin{example}
The curves $$3y^2 = (x^2 - 4) (x^4 + 7x^2 + 1)$$ and
$$- y^2 = (x^2 + 4) (x^4 + 3x^2 + 1)$$ over $\BQ$ are geometrically 
non-isomorphic, and yet their Jacobians are isomorphic to one another
over $\BQ$.
\end{example}

\begin{proof}
If we take the two curves obtained by taking $t = 2$ in 
Theorem~\ref{T-nonsimple2}, replace $x$ by $x/2$ in each equation, and
twist both curves by $2$, we get the two curves given above.
\end{proof}

\begin{example}
The curves
$$5 y^2 = -   6 x^6 -  64 x^5 - 113 x^4 + 262 x^3
          - 331 x^2 + 584 x   + 232 $$
and
$$2 y^2 = -  21 x^6 - 236 x^5 +  45 x^4 - 440 x^3
          - 615 x^2 -  76 x   - 553 $$
are geometrically non-isomorphic, but their Jacobians become isomorphic
to one another over $\BQ(\sqrt{-1})$.  Furthermore, their Jacobians are
absolutely simple.
\end{example}

\begin{proof}
Take $v = 2$ in Theorem~\ref{T-simple2}.  We find that $\rho_1 = -w/4$,
$\rho_2 = w/4$,  $\rho_3 = 2w + 2$, $\rho_4 = w - 1$, 
$\rho_5 = -w - 1$, and $\rho_6 = -2w + 2$, where $w = \sqrt{2}$.  The 
curves $C$ and $C'$ in the theorem are $y^2 = f_1$ and $y^2 = f_2$, 
where
\begin{align*}
f_1 &= -(30625/32) x^6 - (67375/16) x^5 - (305025/64) x^4 \\
    & \qquad  - (23765/16) x^3 + (28665/16) x^2 + (1715/2) x - (735/2)
\end{align*}
and 
$$ f_2 = -(553/2) x^6 + 38 x^5 - (615/2) x^4 + 220 x^3 
         + (45/2) x^2 + 118 x - (21/2).$$
Replacing $x$ with $-2/(x+1)$ in $f_1$ and multiplying the result by 
$(1/5)(2/7)^2 (x+1)^6$ gives rise to the first curve given in the
example.  Replacing $x$ with $-1/x$ in $f_2$ and multiplying the result
by $2 x^6$ gives rise to the second curve.  Thus the Jacobians of the
two curves become isomorphic to one another over
$\BQ(\sqrt{v(v-4)}\,) = \BQ(\sqrt{-1})$. The Jacobians are simple
because we chose our $v$ to be $2$ modulo~$13$.
\end{proof}

\begin{example}
\label{EX-realpoints}
The curves
$$ y^2 + (x^3 + x^2 + x) y =     31 x^6 -  38 x^5 - 217  x^4 - 380 x^3 
                             +  304 x^2 + 501 x   - 366  $$
and
$$11 y^2                   = -   49 x^6 - 378 x^5 - 755  x^4 + 110 x^3
                             - 2285 x^2 + 732 x   - 1368 $$
are geometrically non-isomorphic, but their Jacobians are isomorphic 
to one another over~$\BQ$.  Furthermore, their Jacobians are absolutely
simple.
\end{example}

\begin{proof}
Take $v = -4/3$ in Theorem~\ref{T-simple2}.  The curves $C$ and $C'$
we obtain are $y^2 = f_1$ and $y^2 = f_2$, where
\begin{align*}
f_1 &= (28125/268912) x^6 - (11250/16807) x^5 \\
    & \qquad + (3154875/1882384) x^4 - (812325/470596) x^3 \\
    & \qquad\qquad - (57675/470596) x^2 + (26325/16807) x - (2025/2401)
\end{align*}
and
\begin{align*}
f_2 &= -(131769/38416) x^6 + (11979/343) x^5   \\
    & \qquad - (5595645/38416) x^4 + (62535/196) x^3 \\
    & \qquad\qquad - (3735435/9604) x^2 + (86229/343) x - (23199/343).
\end{align*}
If we replace $x$ with $(2x+2)/(x+2)$ in $f_1$, multiply the result
by $(343/5)^2 (x+2)^6$, and twist by~$3$, we get the curve
$$y^2 =    125 x^6 -  150 x^5 -  865 x^4 - 1518 x^3
        + 1217 x^2 + 2004 x   - 1464;$$
replacing $y$ with $2y + (x^3 + x^2 + x)$ gives the first curve in
the example.  If we replace $x$ with $(x+2)/(x+1)$ in $f_2$, multiply
the result by $(1/11) 196^2 (x+1)^6$, and twist by~$3$, we get the
second curve in the example.   The Jacobians of the two curves become
isomorphic to one another over $\BQ(\sqrt{v(v-4)}\,) = \BQ$.  The
Jacobians are absolutely simple because their reductions modulo $17$
are absolutely simple.
\end{proof}

\begin{rem}
It is easy to see that the first curve in Example~\ref{EX-realpoints}
has real-valued points, while the second curve does not.  It follows
that the real topology of a curve over $\BQ$ is not determined by
its Jacobian.  Furthermore, suppose we choose a positive integer $d$
such that the quadratic twist of the second curve by $d$ has rational
points.  The quadratic twist of the first curve by $d$ will still not
have any real points, let alone any rational points, so we see that the
existence of rational points on a genus-$2$ curve over $\BQ$ is not
determined by its Jacobian, even if the Jacobian is absolutely simple.
\end{rem}

There are also triples $(r,s,t)$ that satisfy the hypotheses of 
Theorem~\ref{T-gensimple} but that do not lie on the elliptic surface 
discussed in Section~\ref{S-application}.

\begin{example}
\label{EX-twiceseen}
The curves
$$y^2 =    x^6          - 24 x^4 + 80 x^3 - 63 x^2 -  24 x -   2$$
and
$$y^2 = -2 x^6 +  6 x^5 +  9 x^4 - 48 x^3          + 162 x - 171$$
are geometrically non-isomorphic, but their Jacobians become isomorphic
to one another over $\BQ(\sqrt{2})$.  Furthermore, their Jacobians are
absolutely simple.
\end{example}

\begin{proof}
We take $r = -7/4$ and $s = 1/2$ and $t = 1/4$ in 
Theorem~\ref{T-gensimple}.  Let $\xi$ be a root of the irreducible
polynomial $$x^6 + 6x^4 + 9x^2 + 16$$ and let $K$
be the number field generated by~$\xi$.  The polynomial $h$ of
Section~\ref{S-Galois} is $$h = T^3 + 3/4 T^2 + 9/16 T + 3/64$$ and its
roots are the elements
\begin{align*}
\beta_1 &= (                 -\xi^4         -7\xi^2        -12)/16 \\
\beta_2 &= (         - \xi^5 +\xi^4 -5\xi^3 +7\xi^2 -10\xi    )/32 \\
\beta_3 &= (\phantom{-}\xi^5 +\xi^4 +5\xi^3 +7\xi^2 +10\xi    )/32
\end{align*}
of $K$.  The polynomials $g_i$ are given by
$$g_i = x^2 - 2\beta_i x + \beta_i^2 / 2 + 3/32,$$ and their roots 
(indexed in accordance with the proof of Theorem~\ref{T-Galois}) are
\begin{align*}
r_1 &=(                 -  \xi^4 -  \xi^3 -  7\xi^2 -  7\xi - 12)/16 \\
r_2 &=(                 -  \xi^4 +  \xi^3 -  7\xi^2 +  7\xi - 12)/16 \\
r_3 &=(         - \xi^5 -  \xi^4 - 3\xi^3 -  3\xi^2 -  8\xi -  8)/32 \\
r_4 &=(         - \xi^5 + 3\xi^4 - 7\xi^3 + 17\xi^2 - 12\xi +  8)/32 \\
r_5 &=(\phantom{-}\xi^5 + 3\xi^4 + 7\xi^3 + 17\xi^2 + 12\xi +  8)/32 \\
r_6 &=(\phantom{-}\xi^5 -  \xi^4 + 3\xi^3 -  3\xi^2 +  8\xi -  8)/32
\end{align*}            
We note that the two subgroups $G$ and $G'$ that appear in the proof of
Theorem~\ref{T-Galois} are indeed Galois stable.  We compute that the 
two Richelot duals are $y^2 = f_1 $ and $y^2 = f_2$, where 
\begin{align*}
f_1 &= - (81/512) x^6 - (1215/1024) x^5 
       - (21141/8192) x^4 - (8991/8192) x^3 \\
    & \qquad  - (19683/131072) x^2 - (2187/262144) x + (729/2097152)\\
f_2 &= - (1863/256) x^6 - (3159/512) x^5 
       - (26973/4096) x^4 - (11421/4096) x^3 \\
    & \qquad - (76545/65536) x^2 - (28431/131072) x - (13851/1048576).
\end{align*}
Evaluating $f_1$ at $(2-x)/(4x)$ and multiplying the result by 
$(256/9)^2 x^6$ gives the first curve mentioned in the example; 
evaluating $f_2$ at $-x/(4x-8)$ and multiplying the result by 
$(128/9)^2 (x-2)^6$ gives the second curve.  These curves are 
geometrically non-isomorphic, and their Jacobians become isomorphic
over $\BQ(\sqrt{st}\,) = \BQ(\sqrt{2})$.  Furthermore, their 
Jacobians are absolutely simple because their reductions modulo $7$
are absolutely simple.
\end{proof}

\begin{rem}
We obtained Example~\ref{EX-twiceseen} from a triple $(r,s,t)$ that
does not lie on the elliptic surface from Section~\ref{S-application},
but the same example can be obtained from the triple 
$(r,s,t) = (-10,-1,-2)$, which \emph{does} lie on the surface.
\end{rem}

We computed all triples $(r,s,t)$ of na\"{\i}ve height at most $20$ for
which the curve $\CC(r,s,t)$ has two $\BQ$-rational Richelot duals
whose Jacobians are absolutely simple and isomorphic over~$\BQ$.  Of
all the examples we found, the triple $(r,s,t) = (-19/3, -6, -1/6)$
gave rise to the curves with the smallest coefficients:

\begin{example}
\label{EX-smallest}
The curves
$$y^2 = -9 x^6 +  6 x^5 -  47 x^4 -  14 x^3 -   5 x^2 - 36 x - 72 $$
and 
$$y^2 =  8 x^6 - 60 x^5 + 235 x^4 - 186 x^3 - 239 x^2 - 30 x -  1 $$
are geometrically non-isomorphic, but their Jacobians are isomorphic 
to one another over~$\BQ$.  Furthermore, their Jacobians are absolutely
simple. \qed
\end{example}
 
\begin{example}
\label{EX-g3}
The Jacobian of the hyperelliptic curve
$$3 v^2 = - 17 u^8 + 56 u^7 - 84 u^6 + 56 u^5 
          - 70 u^4 - 56 u^3 - 84 u^2 - 56 u - 17$$
and the Jacobian of the plane quartic
$$x^4 + 4 y^4 + 4 z^4 + 20 x^2 y^2 - 8 x^2 z^2 + 16 y^2 z^2 = 0$$
are isomorphic to one another over $\BQ$.
\end{example}

\begin{proof}
We take $t=1$ in Theorem~\ref{T-nonsimple3}.  The plane quartic $Q(1)$ 
from the theorem is the plane quartic given in the example.  The 
hyperelliptic curve from the theorem is given by the pair of 
homogeneous equations
\begin{align}
W^2 Z^2 &= -(1/3) X^4 - (4/3) Y^4 + Z^4 \label{E-exH1}\\
      0 &= -X^2 + 2Y^2 + 2Z^2.          \label{E-exH2}
\end{align}
We dehomogenize the equations by setting $Z=1$, and we parametrize
the conic given by Equation~\ref{E-exH2} by setting 
\begin{align*}
 X &= 2(u^2+1)/(u^2 + 2u - 1)\\
 Y &= (u^2 - 2u - 1)/(u^2 + 2u - 1).
\end{align*}
Taking $W = v / (u^2 + 2u - 1)^4$ then gives us the hyperelliptic
curve in our example.
\end{proof}

\begin{rem}
We note that the discriminant of the degree-$8$ polynomial used to 
define the hyperelliptic curve in Example~\ref{EX-g3} is $2^{94}$!
\end{rem}


\end{document}